\documentclass{amsart}
\usepackage{tipa}
\usepackage{amssymb}
\usepackage{stmaryrd}
\usepackage{graphicx}
\usepackage{amsmath,amsfonts}

\newtheorem{theorem}{Theorem}[section]

\theoremstyle{definition}

\newtheorem{remark}[theorem]{Remark}
\newtheorem{question}[theorem]{Question}

\numberwithin{equation}{section}

\begin{document}

\title[Arithmetic properties of polynomials]{Arithmetic properties of polynomials}

\author{Yong Zhang}
\address{School of Mathematics and Statistics, Changsha University of Science and Technology,
Changsha 410114, People's Republic of China}

 \email{zhangyongzju$@$163.com}
\thanks{This research was supported by the National Natural Science Foundation of China (Grant No.~11501052) and the Natural Science Foundation of Zhejiang Province (Grant No.
LQ13A010012).}

\author{Zhongyan Shen}
\address{Department of Mathematics, Zhejiang International Studies University,
Hangzhou, 310012, People's Republic of China }

\email{huanchenszyan$@$yahoo.com.cn}

\subjclass[2010]{Primary 11D25; Secondary 11D72, 11G05}
\date{}

\keywords{Diophantine system, parametric solution, Pell's equation,
elliptic curve}

\begin{abstract} In this paper, first, we prove that the Diophantine
system
\[f(z)=f(x)+f(y)=f(u)-f(v)=f(p)f(q)\]
has infinitely many integer solutions for $f(X)=X(X+a)$ with nonzero
integers $a\equiv 0,1,4\pmod{5}$. Second, we show that the above
Diophantine system has an integer parametric solution for
$f(X)=X(X+a)$ with nonzero integers $a$, if there are integers
$m,n,k$ such that
\[\begin{cases}
\begin{split}
(n^2-m^2)(4mnk(k+a+1)+a(m^2+2mn-n^2))&\equiv0\pmod{(m^2+n^2)^2},\\
(m^2+2mn-n^2)((m^2-2mn-n^2)k(k+a+1)-2amn)&\equiv0\pmod{(m^2+n^2)^2},
\end{split}
\end{cases}\] where $k\equiv0\pmod{4}$ when $a$ is even, and
$k\equiv2\pmod{4}$ when $a$ is odd. Third, we get that the
Diophantine system
\[f(z)=f(x)+f(y)=f(u)-f(v)=f(p)f(q)=\frac{f(r)}{f(s)}\]
has a five-parameter rational solution for $f(X)=X(X+a)$ with
nonzero rational number $a$ and infinitely many nontrivial rational
parametric solutions for $f(X)=X(X+a)(X+b)$ with nonzero integers
$a,b$ and $a\neq b$. At last, we raise some related questions.
\end{abstract}

\maketitle

\section{Introduction}
The $k$th $n$-gonal number is given by
\[P^n_k=\frac{k((n-2)(k-1)+2)}{2}.\]
For $n=3$, $P^3_k=k(k+1)/2$ are called triangular numbers. In 1968,
W. Sierpi\'{n}ski \cite{Sierpinski} showed that there are infinitely
many triangular numbers which at the same time can be written as the
sum, difference, and product of other triangular numbers. For $n=4$,
$P^4_k=k^2$, it is easy to show that $(4m^2 + 1)^2$ is the sum,
difference, and product of squares. In 1986, S. Hirose \cite{Hirose}
proved that for $n=5,6,8$, there are infinitely many $n$-gonal
numbers which at the same time can be written as the sum,
difference, and product of other $n$-gonal numbers. The cases with
$n=7$ and $n\geq9$ are still open.

Some authors proved similar results for the sum and the difference
only, such as R. T. Hansen \cite{Hansen} for $n=5$, W. J.
O$'$Donnell \cite{ODonnell1,ODonnell2} in cases of $n=6,8$, H. J.
Hindin \cite{Hindin} for $n=7$, and S. Ando \cite{Ando} in cases of
general $n$. In 1982, L. C. Eggan, P. C. Eggan and J. L. Selfridge
\cite{Eggan-Eggan-Selfridge} showed that for every $n$, there are
infinitely many $n$-gonal numbers that can be written as the product
of two other $n$-gonal numbers. In 2003, A. F. Beardon
\cite{Beardon} studied the integer solutions $(n,a,b,c,d)$ of the
Diophantine system
\[P(n)=P(a)+P(b)=P(c)-P(d),~~P(a)P(b)P(c)P(d)\neq0,\]
where $P$ is a quadratic polynomial with integer coefficients. Some
related information on $n$-gonal numbers could be found in
\cite[p.1-p.39]{Dickson2}.

A natural question is to consider the integer solutions of the
Diophantine system
\begin{equation}\label{Eq11}
f(z)=f(x)+f(y)=f(u)-f(v)=f(p)f(q),
\end{equation}
where $f(X)$ is a polynomial with rational coefficients and
$\deg{f}\geq2$.

By the theory of Pell's equation, we have:

\begin{theorem}\label{thm11}
Equation (1.1) has infinitely many integer solutions
$(z,x,y,u,v,p,q)$ for $f(X)=X(X+a)$ with nonzero integers $a\equiv
0,1,4\pmod{5}$.
\end{theorem}

\begin{theorem}\label{thm12}
For $f(X)=X(X+a)$ with nonzero integers $a$, if there are integers
$m,n,k$ such that
\[\begin{cases}
\begin{split}
(n^2-m^2)(4mnk(k+a+1)+a(m^2+2mn-n^2))&\equiv0\pmod{(m^2+n^2)^2},\\
(m^2+2mn-n^2)((m^2-2mn-n^2)k(k+a+1)-2amn)&\equiv0\pmod{(m^2+n^2)^2},
\end{split}
\end{cases}\] where $k\equiv0\pmod{4}$ when $a$ is even, and
$k\equiv2\pmod{4}$ when $a$ is odd. Then Equation (1.1) has an
integer parametric solution.
\end{theorem}

To study the integer solutions of (1.1) seems difficult for general
polynomials $f(X)$ with $\deg{f}\geq2$, so we turn to consider the
rational parametric solutions of the Diophantine system
\begin{equation}\label{Eq12}
f(z)=f(x)+f(y)=f(u)-f(v)=f(p)f(q)=\frac{f(r)}{f(s)}.
\end{equation}

For reducible quadratic polynomials, we prove the following
statement.

\begin{theorem}\label{thm13}
For $f(X)=X(X+a)$ with nonzero rational number $a$, Equation (1.2)
has a five-parameter rational solution.
\end{theorem}

For reducible cubic polynomials, by the theory of elliptic curves,
we have:

\begin{theorem}\label{thm14}
 For $f(X)=X(X+a)(X+b)$ with nonzero integers $a,b$ and $a\neq b$, Equation (1.2) has infinitely many rational parametric
solutions.
\end{theorem}

\section{The proofs of Theorems}
\begin{proof}[\textbf{Proof of Theorem 1.1.}] (1)~The cases $a\equiv 0,1\pmod {5}$. Let us start with the equation $f(z) = f(p)f(q)$, where $f(X)=X(X+a)$.
Write $q = p + k$ for some integer $k$, then we obtain
\[(2z+a)^2 = 4(p^2+ap)((p + k)^2 + a(p + k)) + a^2.\]
The polynomial part of the Puiseux expansion of
\[\sqrt{4(p^2+ap)((p + k)^2 + a(p + k)) + a^2}\]
is given by $ak+2(a+k)p+2p^2.$ If there exists a large integral
solution, then
\[2z + a = ak + 2(a + k)p + 2p^2.\]
We get in this
case that $a^2(k+1)(k-1) = 0$, that is $k= \pm1.$ In these cases $q
= p\pm1.$ If $k= 1,$ then we obtain the solutions $z = p^2 + (a+
1)p,$~or $-p^2+(-a-1)p-a.$ If $k=-1,$ then we get $z = p^2 +
(a-1)p-a,$~or $-p^2+(-a+1)p.$

Let us deal with the equation $f(z) = f(x)+f(y)$, where $z =
p^2+(a+1)p$ (We only consider this solution, the other should work
in a similar way). We obtain
\[(2(p^2 + (a + 1)p) + a)^2 = (2x + a)^2 + (2y + a)^2-a^2.\]
Take $y = x + b$ for some integer $b$. It follows that
\[(2(p^2 + (a + 1)p) + a)^2-2(2x + a + b)^2 = 2b^2-a^2.\]
Let $X=2(p^2 + (a + 1)p) + a,Y=2x + a + b$, we get the Pell's
equation
\begin{equation}\label{Eq21}
X^2-2Y^2=2b^2-a^2.
\end{equation}
It is easy to provide infinitely many parameterized solutions by the
formula
\begin{equation}\label{Eq22}
X + Y\sqrt{2} =(1+\sqrt{2})^{2m+1}(a+b\sqrt{2}),~m\in
\mathbb{Z}.\end{equation}

When $m = 0$, it yields the trivial solution with $x = 0$. When $m =
1$, we get
\[2x+a+b=5a+7b,2(p^2 +(a + 1)p)+a=7a+10b,\]
then $x = 2a + 3b$ and
\[(2p + a + 1)^2-(a+7)^2=20b-48.\]
From the second formula, we obtain
\[(p-3)(p+a+4)=5b-12.\]
If we write $5b-12 = b_1b_2$, then we have $p = b_1 +3$ and $b_2 =
b_1 +a +7$ as a solution. Here we need that
\[b = \frac{b_1(b_1+ a+ 7) + 12}{5}\] is an integer, it is if $a\equiv 0\pmod {5},~b_1\equiv 1,2 \pmod {5}$, or $a\equiv 1 \pmod {5},~b_1\equiv
3,4 \pmod {5}$.

Up to now we constructed infinitely many integer solutions of the
equations \[f(z) = f(x) + f(y) = f(p)f(q),\] so it remains to
consider the case $f(z) = f(u)-f(v)$, where $z =
(b_1+3)^2+(a+1)(b_1+3)$. Let $v = u-c$ for some integer $c$. Then
$f(z) =c(2u+a-c)$, a linear equation in $u$. Hence,
\[u=\frac{z^2+az+c^2-ac}{2c}.\] As a solution, fix $b_1 = 2c-3$ and
we obtain
\[\begin{split}
&u = (4c^2+2ac+a+4c)(2c+a)+\frac{a+5c}{2},\\
&v = (4c^2+2ac+a+4c)(2c+a)+\frac{a+3c}{2},\\
&z = 2c(2c+a+1).
\end{split}\]

According to the other variables we have
\[\begin{split}
&p = 2c,\\
&q = 2c + 1,\\
&x =2a+3\frac{4c^2+2ac-3a+2c}{5},\\
&y =2a+4\frac{4c^2+2ac-3a+2c}{5}.
\end{split}\]

To get integral values of $x,y,u,v$, we need the following
conditions
\[a\equiv c \pmod{2}\] and
\[4c^2+2ac-3a+2c\equiv0\pmod {5}.\]

(i) If $a\equiv 0\pmod {5},~b_1\equiv 1,2 \pmod {5}$, from the
second condition, we have
\[c\equiv0,2\pmod {5}.\] From the first condition, we get
\[a\equiv 0\pmod {10},~c\equiv0,2\pmod {10},\]
and \[a\equiv 5\pmod {10},~c\equiv5,7\pmod {10}.\] For $a\equiv
0\pmod {10}$, take $c=10t,$~or~$10t+2$, where $t$ is an integer
parameter, we have
\[b=80t^2+4t+4at-\frac{3a}{5},~or~80t^2+36t+4+4at+\frac{a}{5}\in
\mathbb{Z}.\]  For \mbox{$a\equiv 5\pmod {10}$}, put
$c=10t+5,$~or~$10t+7$, we get
\[b=80t^2+84t+22+4at+\frac{7a}{5},~or~80t^2+116t+42+4at+\frac{11a}{5}\in \mathbb{Z}.\]

(ii) If $a\equiv 1\pmod {5},~b_1\equiv 3,4 \pmod {5}$, from the
second condition, we have
\[c\equiv1,3\pmod {5}.\] From the first condition, we obtain
\[a\equiv 1\pmod {10},~c\equiv1,3\pmod {10},\]
and\[a\equiv 6\pmod {10},~c\equiv6,8\pmod {10}.\] For $a\equiv
1\pmod {10}$, take $c=10t+1,$~or~$10t+3$, where $t$ is an integer
parameter, we have
\[b=80t^2+20t+1+4at+\frac{1-a}{5},~or~80t^2+52t+8+4at+\frac{2+3a}{5}\in \mathbb{Z}.\]
For $a\equiv 6\pmod {10}$, put $c=10t+6,$~or~$10t+8$, we get
\[b=80t^2+100t+31+4at+\frac{1+9a}{5},~or~80t^2+132t+54+4at+\frac{2+13a}{5}\in \mathbb{Z}.\]

(2)~The case $a\equiv 4\pmod {5}$. Let us note that Equation (2.1)
has another family of integer solutions
\[X + Y\sqrt{2} =(1+\sqrt{2})^{2m+1}(-a+b\sqrt{2}),~m\in \mathbb{Z}.\]
When $m = 1$, we get
\[2x+a+b=-5a+7b,2(p^2+(a+1)p)+a=-7a+10b,\]
then $x =-3a + 3b$ and
\[(2p+a+1)^2-(a-7)^2=20b-48.\]
From the second formula, we obtain
\[(p+4)(p+a-3)=5b-12.\]
Writing $5b-12 = b_1b_2$, then we have $p=b_1-4$ and $b_2=b_1+a -7$
as a solution. Here we need that
\[b = \frac{b_1(b_1+a-7) + 12}{5}\] is an integer, it is if $a\equiv 0\pmod {5},~b_1\equiv 3,4 \pmod {5}$, or $a\equiv 4\pmod {5},~b_1\equiv
1,2 \pmod {5}$. We only consider the case $a\equiv 4\pmod {5}$ in
the following.

As in the first part, take $b_1=2c+3$ in $z =
(b_1-4)^2+(a+1)(b_1-4)$, we obtain
\[\begin{split}
&z = (2c-1)(2c+a),\\
&x =-3a+3\frac{4c^2+2ac+3a-2c}{5},\\
&y =-3a+4\frac{4c^2+2ac+3a-2c}{5},\\
&u = (4c^2+2ac-a-4c)(2c+a)+\frac{a+5c}{2},\\
&v = (4c^2+2ac-a-4c)(2c+a)+\frac{a+3c}{2},\\
&p = 2c-1,\\
&q = 2c.
\end{split}\]

To get integral values of $x,y,u,v$, we need the following
conditions
\[a\equiv c \pmod{2}\] and
\[4c^2+2ac+3a-2c\equiv0\pmod {5}.\]

If $a\equiv 4\pmod {5},~b_1\equiv 1,2 \pmod {5}$, from the second
condition, we have
\[c\equiv2,4\pmod {5}.\] From the first condition, we get
\[a\equiv 4\pmod {10},~c\equiv2,4\pmod {10},\]
and\[a\equiv 9\pmod {10},~c\equiv7,9\pmod {10}.\]
As in the above, we
can take $c=10t+2,$~or~$10t+4$ for $a\equiv 4\pmod {10}$, and
$c=10t+7,$~or~$10t+9$ for $a\equiv 9\pmod {10}$, where $t$ is an
integer parameter.

Combining (1) and (2) completes the proof of Theorem 1.1.
\end{proof}

\textbf{Example 1.} When $a = 1$, $c=10t+1$, then Equation (1.1) has
solutions
\[\begin{split}
&z = 400t^2 + 120t+ 8,\\
&x = 240t^2 + 72t + 5,\\
&y = 320t^2 + 96t + 6,\\
&u =8000t^3+4000t^2+665t+36,\\
&v =8000t^3+4000t^2+655t+35,\\
&p = 20t + 2,\\
&q = 20t + 3,
\end{split}\]
where $t$ is an integer parameter.

\begin{remark}
It's worth to note that we have other possibilities to obtain
integer solutions if $a\not \equiv 0,1,4\pmod {5}$. When we apply
the same idea using the other family of solutions with $z
=-p^2+(-a-1)p-a,~p^2 + (a-1)p-a,$~or~$-p^2+(-a+1)p$ for $m=1,$ we do
not obtain new cases of solutions.

Also we only used the first non-trivial value of $m$ in case of the
general Pell's equation. To cover more classes one has to go in
these directions. For example, when $m=2,$ from formula (2.2), we
have
\[2x+a+b=29a+41b,~2(p^2+(a+1)p)+a=41a+58b,\]
then $x=14a+20b$ and
\[(2p+a+1)^2-(a+41)^2 = 116b-1680.\]
Hence, we get \[(p-20)(p+a+21)=29b-420.\] Let $29b-420=b_1b_2,$ then
we have $p=b_1+20,b_2=b_1+a+41.$ Hence,
\[b=\frac{b_1(b_1+a+41)+420}{29}.\]
To make $b$ be an integer, we have
\[a \equiv 0, 2, 3, 5, 6, 8, 10, 11, 15, 19, 23, 24, 26, 28\pmod{29}.\]

As in Theorem 1.1, fix $b_1 = 2c-20$ in $z =
(b_1+20)^2+(a+1)(b_1+20)$ and we obtain
\[\begin{split}
&z = 2c(2c+a+1),\\
&x = 14a+20\frac{4c^2+2ac-20a+2c}{29},\\
&y = 14a+21\frac{4c^2+2ac-20a+2c}{29},\\
&u = (4c^2+2ac+a+4c)(2c+a)+\frac{a+5c}{2},\\
&v = (4c^2+2ac+a+4c)(2c+a)+\frac{a+3c}{2},\\
&p = 2c,\\
&q = 2c+1.
\end{split}\]

To get integral values of $x,y,u,v$, we need the following
conditions
\[a\equiv c \pmod{2}\] and
\[4c^2+2ac-20a+2c\equiv0\pmod {29}.\]

In the following table, we give the conditions of $a,b_1,c$ such
that $b$ is an integer.
\[\begin{tabular}{|c|c|c|}
\hline

 $a$ & $b_1$ & $c$ \\

\hline
$\equiv 0\pmod{58}$  & $\equiv 8,9\pmod{29}$ & $\equiv0,14\pmod{58}$ \\
\hline
 $\equiv 29\pmod{58}$  & $\equiv 8,9\pmod{29}$ & $\equiv 29,43\pmod{58}$ \\

\hline
$\equiv 2\pmod{58}$ & $\equiv 1,14\pmod{29}$ & $\equiv 46,54\pmod{58}$ \\
\hline
$\equiv 31\pmod{58}$ & $\equiv 1,14\pmod{29}$ & $\equiv 17,25\pmod{58}$ \\

\hline
$\equiv 3\pmod{58}$ & $\equiv 15,28\pmod{29}$ & $\equiv 3,53\pmod{58}$ \\
\hline
$\equiv 32\pmod{58}$ & $\equiv 15,28\pmod{29}$ & $\equiv 32,24\pmod{58}$ \\

\hline
$\equiv 5\pmod{58}$ & $\equiv 20,21\pmod{29}$  & $\equiv  35,49\pmod{58}$\\
\hline
$\equiv 34\pmod{58}$ & $\equiv 20,21\pmod{29}$  & $\equiv  6,20\pmod{58}$\\

\hline
$\equiv 6\pmod{58}$ & $\equiv 17,23\pmod{29}$ & $\equiv 4,36\pmod{58}$\\
\hline
$\equiv 35\pmod{58}$ & $\equiv 17,23\pmod{29}$ & $\equiv 33,7\pmod{58}$\\

\hline
$\equiv 8\pmod{58}$ & $\equiv 2,7\pmod{29}$ & $\equiv  40,28\pmod{58}$\\
\hline
$\equiv 37\pmod{58}$ & $\equiv 2,7\pmod{29}$ & $\equiv  11,57\pmod{58}$\\

\hline
$\equiv 10\pmod{58}$ & $\equiv 11,25\pmod{29}$ & $\equiv  30,8\pmod{58}$\\
\hline
$\equiv 39\pmod{58}$ & $\equiv 11,25\pmod{29}$ & $\equiv  1,37\pmod{58}$\\

\hline
$\equiv 11\pmod{58}$ & $\equiv 16,19\pmod{29}$ & $\equiv  5,47\pmod{58}$\\
\hline
$\equiv 40\pmod{58}$ & $\equiv 16,19\pmod{29}$ & $\equiv  34,18\pmod{58}$\\

\hline
$\equiv 15\pmod{58}$ & $\equiv 5,26\pmod{29}$ & $\equiv  23,27\pmod{58}$\\
\hline
$\equiv 44\pmod{58}$ & $\equiv 5,26\pmod{29}$ & $\equiv  52,56\pmod{58}$\\

\hline
$\equiv 19\pmod{58}$ & $\equiv 3,24\pmod{29}$ & $\equiv  51,55\pmod{58}$\\
\hline
$\equiv 48\pmod{58}$ & $\equiv 3,24\pmod{29}$ & $\equiv  22,26\pmod{58}$\\

\hline
$\equiv 23\pmod{58}$ & $\equiv 10,13\pmod{29}$ & $\equiv 31,15\pmod{58}$\\
\hline
$\equiv 52\pmod{58}$ & $\equiv 10,13\pmod{29}$ & $\equiv 2,44\pmod{58}$\\

\hline
$\equiv 24\pmod{58}$ & $\equiv 4,18\pmod{29}$ & $\equiv  12,48\pmod{58}$\\
\hline
$\equiv 53\pmod{58}$ & $\equiv 4,18\pmod{29}$ & $\equiv  41,19\pmod{58}$\\

\hline
$\equiv 26\pmod{58}$ & $\equiv 22,27\pmod{29}$ & $\equiv  38,50\pmod{58}$\\
\hline
$\equiv 55\pmod{58}$ & $\equiv 22,27\pmod{29}$ & $\equiv  9,21\pmod{58}$\\

\hline
$\equiv 28\pmod{58}$ & $\equiv 6,12\pmod{29}$ & $\equiv 42,16\pmod{58}$\\
\hline
$\equiv 57\pmod{58}$ & $\equiv 6,12\pmod{29}$ & $\equiv 13,45\pmod{58}$\\

\hline

\end{tabular}
\]
For example, when $a\equiv 0\pmod{58},$ we can take
$c=58t,$~or~$58t+14$, then
\[b=464t^2+4t+4at-\frac{20a}{29},~or~464t^2+228t+28+4at+\frac{8a}{29} \in \mathbb{Z}.\]

In view of Theorem 1.1, the cases $a\equiv 2,3\pmod{5}$ are not
covered. To get more $a$ such that Equation (1.1) has infinitely
many integer solutions for $f(X)=X(X+a)$, we use the congruences
modulo 29 in the above table. Solving the congruences
\[a\equiv 2\pmod{5},a\equiv 0, 2,
3, 5, 6, 8, 10, 11, 15, 19, 23, 24, 26, 28\pmod{29},\] we have
\[a\equiv 87, 2, 32, 92, 122, 37, 97, 127, 102, 77, 52, 82, 142, 57 \pmod{145}.\]
By the congruences
\[a\equiv 3\pmod{5},a\equiv 0, 2, 3, 5, 6, 8, 10,
11, 15, 19, 23, 24, 26, 28\pmod{29},\] we get
\[a\equiv 58, 118, 3, 63, 93, 8, 68, 98, 73, 48, 23, 53, 113, 28 \pmod{145}.\]
Then if
\[\begin{split}
a\equiv &2, 3, 8, 23, 28, 32, 37, 48, 52, 53, 57, 58, 63, 68, 73,\\
&77, 82, 87, 92, 93, 97, 98, 102, 113, 118, 122, 127, 142
\pmod{145},
\end{split}\] (1.1) has infinitely many integer
solutions for $f(X)=X(X+a)$.

Similar congruences can be obtained for other values of $m$, and
combining these systems via the Chinese remainder theorem would
cover almost all classes. However, it seems difficult to cover all
integers $a\equiv 2,3\pmod{5}$.
\end{remark}

Note that Example 1 gives an integer parametric solution of Equation
(1.1) with $f(X)=X(X+1)$, we try to generalize these formulas and
get the proof of Theorem 1.2.

\begin{proof}[\textbf{Proof of Theorem 1.2.}]
Firstly, we study the equation $f(z)=f(p)f(q)$ for $f(X)=X(X+a).$
Take $p=At+k,q=At+k+1,$ then we have
\[z=(At+k)(At+k+a+1).\]

Secondly, we consider the equation $f(z)=f(u)-f(v).$ Let
$u=Bt^3+Ct^2+Dt+E,v=Bt^3+Ct^2+Ft+G$, then
\[\begin{split}
(At+k)(At+k+a)&(At+k+1+a)(At+k+1)\\
&=(Dt-Ft+E-G)(2Bt^3+2Ct^2+Dt+Ft+E+G+a).
\end{split}\]
To determinate the coefficients of $u,v$, by the method of
undetermined coefficients, we obtain
\begin{equation}
\begin{split}
A^4&=2B(D-F),\\
2A^3(a+2k+1)&=2BE-2BG+2CD-2CF,\\
A^2(a^2+6ak+6k^2+3a+6k+1)&=2CE-2CG+D^2-F^2, \\
A(a+2k+1)(2ak+2k^2+a+2k)&=2DE+Da-2FG-Fa, \\
k(k+1)(a+k+1)(a+k)&=(E-G)(E+a+G).
\end{split}
\end{equation}
In order to find a solution of Equation (2.3), let
$B=A^3,E-G=\frac{k}{2},$ then \[F=D-\frac{A}{2},G=E-\frac{k}{2}.\]
Put $B,F,G$ into Equation (2.3), and solve it for $C,D,E,$ we get
\[\begin{split}
C&=A^2(2a+3k+2),\\
D&=A(a^2+4ak+3k^2+3a+4k)+\frac{5A}{4},\\
E&=(ak+k^2+a+2k)(a+k)+\frac{2a+5k}{4}.
\end{split}\]
Hence,
\[\begin{split}
u=&A^3t^3+A^2(2a+3k+2)t^2+\bigg(A(a^2+4ak+3k^2+3a+4k)+\frac{5A}{4}\bigg)t\\
&+(ak+k^2+a+2k)(a+k)+\frac{2a+5k}{4},\\
v=&A^3t^3+A^2(2a+3k+2)t^2+\bigg(A(a^2+4ak+3k^2+3a+4k)+\frac{3A}{4}\bigg)t\\
&+(ak+k^2+a+2k)(a+k)+\frac{2a+3k}{4}.
\end{split}\]

At last, we study the equation $f(z)=f(x)+f(y).$ Put
$x=Ht^2+It+J,y=Kt^2+Lt+M$, then
\[\begin{split}
(At+k)&(At+k+a)(At+k+1+a)(At+k+1)\\
=&(H^2+K^2)t^4+(2IH+2KL)t^3+(2HJ+Ha+2KM+Ka+L^2-I^2)t^2\\
&+(2IJ+2LM+La+Ia)t+J^2+Ja+M^2+Ma.
\end{split}\]
To determinate the coefficients of $x,y$, by the method of
undetermined coefficients, we obtain
\begin{equation}
\begin{split}
A^4&=H^2+K^2,\\
2A^3(a+2k+1)&=2IH+2KL,\\
A^2(a^2+6ak+6k^2+3a+6k+1)&=2HJ+Ha+2KM+Ka+L^2-I^2, \\
A(a+2k+1)(2ak+2k^2+a+2k)&=2IJ+2LM+La+Ia, \\
k(k+1)(a+k+1)(a+k)&=J^2+Ja+M^2+Ma.
\end{split}
\end{equation}
Solve the first equation of Equation (2.4), we get an integer
parametric solution
\[A=n^2+m^2,H=4nm(n^2-m^2),K=m^4-6m^2n^2+n^4,\]
where $n>m$ are integer parameters. Take $A,H,K$ into the second,
third and fourth equations of Equation (2.4), and solve them for
$I,L,M,$ then
\[\begin{split}
I&=\frac{4nm(n^2-m^2)(a+2k+1)}{n^2+m^2},\\
L&=\frac{(m^4-6m^2n^2+n^4)(a+2k+1)}{n^2+m^2},\\
M&=\frac{4nm(n^2-m^2)J+(m^2+n^2)^2k^2+(m^2+n^2)^2(a+1)k+2amn(m^2+2mn-n^2)}{m^4-6m^2n^2+n^4}.
\end{split}\]
Put $M$ into the firth equation of Equation (2.4), we have
\[\begin{split}
J=\frac{(n^2-m^2)(4mnk(k+a+1)+a(m^2+2mn-n^2))}{(n^2+m^2)^2}.
\end{split}\]
Hence,
\[\begin{split}
M=\frac{(m^2+2mn-n^2)((m^2-2mn-n^2)k(k+a+1)-2amn)}{(n^2+m^2)^2}.
\end{split}\]
So
\[\begin{split}
x=&4nm(n^2-m^2)t^2+\frac{4nm(n^2-m^2)(a+2k+1)}{n^2+m^2}t\\
&+\frac{(n^2-m^2)(4mnk(k+a+1)+a(m^2+2mn-n^2))}{(n^2+m^2)^2},\\
y=&(m^4-6m^2n^2+n^4)t^2+\frac{(m^4-6m^2n^2+n^4)(a+2k+1)}{n^2+m^2}t\\
&+\frac{(m^2+2mn-n^2)((m^2-2mn-n^2)k(k+a+1)-2amn)}{(n^2+m^2)^2}.
\end{split}\]
According to the other variables, we have
\[\begin{split}
p=&(n^2+m^2)t+k,\\
q=&(n^2+m^2)t+k+1,\\
u=&(n^2+m^2)^3t^3+(n^2+m^2)^2(2a+3k+2)t^2+\bigg((n^2+m^2)(a^2+4ak\\
&+3k^2+3a+4k)+\frac{5(n^2+m^2)}{4}\bigg)t+(ak+k^2+a+2k)(a+k)+\frac{2a+5k}{4},\\
v=&(n^2+m^2)^3t^3+(n^2+m^2)^2(2a+3k+2)t^2+\bigg((n^2+m^2)(a^2+4ak\\
&+3k^2+3a+4k)+\frac{3(n^2+m^2)}{4}\bigg)t+(ak+k^2+a+2k)(a+k)+\frac{2a+3k}{4}.
\end{split}\]
To get integral values of $x,y,u,v$, we can take $t=4(n^2+m^2)T$,
where $T$ is an integer parameter, $k\equiv0\pmod{4}$ when $a$ is
even, and $k\equiv2\pmod{4}$ when $a$ is odd and the following
congruence conditions
\[
\begin{cases}
\begin{split}
(n^2-m^2)(4mnk(k+a+1)+a(m^2+2mn-n^2))&\equiv0\pmod{(m^2+n^2)^2},\\
(m^2+2mn-n^2)((m^2-2mn-n^2)k(k+a+1)-2amn)&\equiv0\pmod{(m^2+n^2)^2}.
\end{split}
\end{cases}\] This completes the proof of Theorem 1.2. \end{proof}

\textbf{Example 2.} (1) When $a=1$, $f(X)=X(X+1)$, we take
$(m,n)=(1,2),$ $t=20T,k=4k_1+2$, then we get the conditions

\[\begin{cases}\begin{split}
384k_1^2+576k_1+195&\equiv0\pmod{25},\\
-112k_1^2-168k_1-60&\equiv0\pmod{25}.
\end{split}
\end{cases}\] Solve these two congruences, we obtain
\[k_1\equiv5,6\pmod{25}.\]
Then
\[k\equiv22,26\pmod{100}.\]
If we put $k=100S+22,$ where $S$ is an integer parameter, then there
are infinitely many integer parametric solutions
\[\begin{split}
z=&8(50T+50S+11)(25T+6+25S),\\
x=&9600T^2+(19200S+4416)T+9600S^2+4416S+507,\\
y=&-2800T^2-(5600S+1288)T-2800S^2-1288S-148,\\
u=&1000000T^3+(3000000S+700000)T^2+(3000000S^2+1400000S\\
&+163325)T+1000000S^3+700000S^2+163325S+12701,\\
v=&5(200T+47+200S)(40T+9+40S)(25T+6+25S),\\
p=&100T+100S+22,\\
q=&100T+100S+23.
\end{split}\]

(2) When $a=2$, $f(X)=X(X+2)$, we take $(m,n)=(1,4),$
$t=68T,k=4k_1$, then we have the conditions
\[\begin{cases}
\begin{split}
3840k_1^2+2880k_1-210&\equiv0\pmod{289},\\
2576k_1^2+1932k_1+112&\equiv0\pmod{289}.
\end{split}
\end{cases}\]
Solve these two congruences, we obtain
\[k_1\equiv104,112\pmod{289}.\]
Then
\[k\equiv416,448\pmod{1156}.\]
If we set $k=1156S+416,$ where $S$ is an integer parameter, then
there exist infinitely many integer parametric solutions
\[\begin{split}
z=&4(289T+289S+104)(1156T+419+1156S),\\
x=&1109760T^2+(2219520S+801600)T+1109760S^2+801600S+144750,\\
y=&744464T^2+(1488928S+537740)T+744464S^2+537740S+97104,\\
u=&1544804416T^3+(4634413248S+1675765344)T^2+(4634413248S^2\\
&+3351530688S+605941965)T+1544804416S^3+1675765344S^2\\
&+605941965S+73034317,\\
v=&(2312T+837+2312S)(668168S^2+1336336ST+668168T^2+482919S\\
&+482919T+87257),\\
p=&1156T+1156S+416,\\
q=&1156T+1156S+417.
\end{split}\]

(3) When $a=3$, $f(X)=X(X+3)$, we take $(m,n)=(1,5),$
$t=104T,k=4k_1+2$, then we get the conditions
\[\begin{cases}\begin{split}
7680k_1^2+15360k_1+4752&\equiv0\pmod{676},\\
7616k_1^2+15232k_1+6132&\equiv0\pmod{676}.
\end{split}
\end{cases}\]
Solve these two congruences, we obtain
\[k_1\equiv 46, 121, 215, 290, 384, 459, 553, 628\pmod{676}.\]
Then
\[k\equiv 186, 486, 862, 1162, 1538, 1838, 2214, 2514\pmod{2704}.\]
If we put $k=2704S+186,$ where $S$ is an integer parameter, then
there are infinitely many integer parametric solutions
\[\begin{split}
z=&4(1352T+1352S+93)(1352T+95+1352S),\\
x=&5191680T^2+(10383360S+721920)T+5191680S^2+721920S+25092,\\
y=&5148416T^2+(10296832S+715904)T+5148416S^2+715904S+24885,\\
u=&19770609664T^3+(59311828992S+4138374656)T^2+(59311828992S^2\\
&+8276749312S+288741908)T+19770609664S^3+4138374656S^2\\
&+288741908S+6715215,\\
v=&2(5408T+379+5408S)(1827904S^2+3655808ST+1827904T^2+254514S\\
&+254514T+8859),\\
p=&2704T+2704S+186,\\
q=&2704T+2704S+187.
\end{split}\]

\begin{proof}[\textbf{Proof of Theorem 1.3.}] For $f(X)=X(X+a)$, let $z=T$, the first equation of Equation (1.2) reduces to
\[T(T+a)=x(x+a)+y(y+a).\]
This can be parameterized by
\[x=-\frac{2Tk+ak+a}{k^2+1},y=-\frac{(k+1)(Tk+ak-T)}{k^2+1},\]
where $k$ is a rational parameter.

From $T(T+a)=u(u+a)-v(v+a),$ we get
\[u=-\frac{Tt^2+at^2-at+T}{t^2-1},v=-\frac{2Tt+at-a}{t^2-1},\]
where $t$ is a rational parameter.

For $T(T+a)=p(p+a)q(q+a)$, put $p=wT$, then
\[T=-\frac{a(aqw+q^2w-1)}{aqw^2+q^2w^2-1},p=-\frac{aw(aqw+q^2w-1)}{aqw^2+q^2w^2-1},\]
where $w$ is a rational parameter.

Take $s=mr,$ from \[T(T+a)=\frac{r(r+a)}{s(s+a)},\] we obtain
\[r=-\frac{a(T^2m+Tam-1)}{T^2m^2+Tam^2-1},s=-\frac{am(T^2m+Tam-1)}{T^2m^2+Tam^2-1},\]
where $m$ is a rational parameter.

Put \[T=-\frac{a(aqw+q^2w-1)}{aqw^2+q^2w^2-1}\] into $x,y,u,v,r,s,$
then Equation (1.2) has a five-parameter rational solution.
\end{proof}

Now we provide the proof of Theorem 1.4.

\begin{proof}[\textbf{Proof of Theorem 1.4.}] To prove this theorem, we need to consider four Diophantine equations. The first one is
\begin{equation}\label{Eq23}
z(z+a)(z+b)=x(x+a)(x+b)+y(y+a)(y+b).
 \end{equation}
Let $z=T,$ and consider Equation (2.5) as a cubic curve with
variables $x,y$:
\[C_1:~x(x+a)(x+b)+y(y+a)(y+b)-T(T+a)(T+b)=0.\]
By the method described in \cite[p.477]{Cohen1}, using Magma
\cite{Bosma-Cannon-Playoust}, $C_1$ is birationally equivalent to
the elliptic curve
\[\begin{split}
E_1:~Y^2=&X^3-432(a^2-ab+b^2)^2X-314928T^6-629856(a+b)T^5\\
&-314928(a^2+4ab+b^2)T^4+23328(a+b)(4a^2-37ab+4b^2)T^3\\
&+1164(8a^4-4a^3b-51a^2b^2-4ab^3+8b^4)T^2+46656ab(2a-b)(a-2b)\\
&\times
(a+b)T-1728(a^2-4ab+b^2)(a^2+2ab-2b^2)(2a^2-2ab-b^2).\end{split}\]
The map $\varphi_1: C_1\rightarrow E_1$ is
\[\begin{split}
X=&\frac{12(3(a^2-ab+b^2)(x+y)+27T(T+a)(T+b)-2(a+b)(a^2-4ab+b^2))}{2a+2b+3x+3y},\\
Y=&\frac{-108(3T+2a-b)(3T+2a+2b)(3T-a+2b)(x-y)}{2a+2b+3x+3y},
\end{split}\]
and its inverse map $\varphi_1^{-1}: E_1\rightarrow C_1$ is
\[\begin{split}
x=&\frac{-(6a+6b)X-Y+972T(T+a)(T+b)-72(a+b)(a^2-4ab+b^2)}{18(X-12a^2+12ab-12b^2)},\\
y=&\frac{-(6a+6b)X+Y+972T(T+a)(T+b)-72(a+b)(a^2-4ab+b^2)}{18(X-12a^2+12ab-12b^2)}.
\end{split}\]
The discriminant of $E_1$ is nonzero as an element of
$\mathbb{Q}(T)$, then $E_1$ is smooth.

Note that the point $(x,y)=(0,T)$ lies on $C_1$, by the map
$\varphi_1$, the corresponding point on $E_1$ is
\[W=(108T^2+36(a+b)T-12a^2+48ab-12b^2,108(3T-b+2a)(3T+2b-a)T).\]
By the group law, we have
\[\begin{split}
[2]W=&\bigg(12(9T^4+12(a+b)T^3+(5a^2+16ab+5b^2)T^2+6ab(a+b)T+3a^2b^2)/T^2,\\
&-108(9T^6+21(a+b)T^5+(16a^2+41ab+16b^2)T^4\\
&+2(a+b)(2a^2+13ab+2b^2)T^3+2ab(4a^2+11ab+4b^2)T^2\\&+6a^2b^2(a+b)T+2a^3b^3)/T^3\bigg).
\end{split}\]

A quick computation reveals that the remainder of the division of
the numerator by the denominator of the $X$-th coordinate of $[2]W$
with respect to $T$ is equal to\[(72a^2b+72ab^2)T+36a^2b^2\] and
thus is nonzero as an element of $\mathbb{Q}(T)$ provided $ab\neq0$.
By a generalization of Nagell-Lutz theorem (see
\cite[p.268]{Connell}), $[2]W$ is of infinite order on $E_1$, then
there are infinitely many $\mathbb{Q}(T)$-rational points on $E_1$.

For $m=2,3,...,$ compute the points $[m]W$ on $E_1$, next calculate
the corresponding point $\varphi_1^{-1}([m]W) = (x_m, y_m)$ on
$C_1$. Then we get infinitely many $\mathbb{Q}(T)$-rational
solutions $(x,y)$ of Equation (2.5).

The second one is
\begin{equation}\label{Eq24}
z(z+a)(z+b)=u(u+a)(u+b)-v(v+a)(v+b).
 \end{equation}
Take $z=T,$ and consider Equation (2.6) as a cubic curve with
variables $u,v$:
\[C_2:~u(u+a)(u+b)-v(v+a)(v+b)-T(T+a)(T+b)=0.\]
By the method described in \cite[p.477]{Cohen1}, using Magma
\cite{Bosma-Cannon-Playoust}, $C_2$ is birationally equivalent with
the elliptic curve
\[\begin{split}
E_2:~V^2=&U^3-432(a^2-ab+b^2)^2U-314928T^6-629856(a+b)T^5\\
&-314928(a^2+4ab+b^2)T^4-629856ab(a+b)T^3-314928a^2b^2T^2\\
&+3456(a^2-ab+b^2)^3.\end{split}\] The map $\varphi_2:
C_2\rightarrow E_2$ is
\[\begin{split}
U=&\frac{12((a^2-ab+b^2)(u-v)+9T(T+a)(T+b))}{u-v},\\
V=&\frac{324T(T+a)(T+b)(2a+2b+3u+3v)}{u-v},
\end{split}\]
and its inverse map $\varphi_2^{-1}: E_2\rightarrow C_2$ is
\[\begin{split}
u=&\frac{-(6a+6b)U-V+972T(T+a)(T+b)+72(a^3+b^3)}{18(U-12a^2+12ab-12b^2)},\\
v=&\frac{-(6a+6b)U+V+972T(T+a)(T+b)+72(a^3+b^3)}{-18(U-12a^2+12ab-12b^2)}.
\end{split}\]
The discriminant of $E_2$ is nonzero as an element of
$\mathbb{Q}(T)$, then $E_2$ is smooth.

Note that the point $(u,v)=(T,0)$ lies on $C_2$, by the map
$\varphi_2$, the corresponding point on $E_2$ is
\[W'=(108T^2+108(a+b)T+12a^2+96ab+12b^2,324(T+a)(T+b)(3T+2a+2b)).\]
By the group law, we get
\[\begin{split}
[2]W'=&\bigg(12(81T^4+108(a+b)T^3+(45a^2+144ab+45b^2)T^2+(12a^3+54a^2b\\
&+54ab^2+12b^3)T+4a^4+4a^3b+27a^2b^2+4ab^3+4b^4)/(3T+2a+2b)^2,\\
&-324(81T^6+135(a+b)T^5+(54a^2+189ab+54b^2)T^4-(12a^3-18a^2b\\
&-18ab^2+12b^3)T^3-(8a^4+68a^3b+66a^2b^2+68ab^3+8b^4)T^2-6ab\\
&\times(a+b)(4a^2+5ab+4b^2)T-6a^2b^2(2a^2+ab+2b^2))/(3T+2a+2b)^3\bigg).
\end{split}\]

By the same method of above, we have infinitely many
$\mathbb{Q}(T)$-rational solutions $(u,v)$ of Equation (2.6).

The third one is
\begin{equation}\label{Eq25}
z(z+a)(z+b)=p(p+a)(p+b)q(q+a)(q+b).
 \end{equation}
Put $z=T,q(q+a)(q+b)=Q$ and consider Equation (2.7) as a cubic curve
with variables $T,p$:
\[C_3:~T(T+a)(T+b)-Qp(p+a)(p+b)=0.\]
By the method described in \cite[p.477]{Cohen1},  using Magma
\cite{Bosma-Cannon-Playoust}, $C_3$ is birationally equivalent to
the elliptic curve
\[\begin{split}
E_3:~Y^2=&X^3-432(a^2-ab+b^2)^2Q^2X+27a^2b^2(a-b)^2Q^2\\
&2(2a-b)^2(a-2b)^2(a+b)^2Q+27a^2b^2(a-b)^2.\end{split}\] Because the
map $\varphi_3: C_3\rightarrow E_3$ is complicated, we omit it. The
discriminant of $E_3$ is nonzero as an element of $\mathbb{Q}(Q)$,
then $E_3$ is smooth.

Note that the point $(T,p)=(a,a)$ lies on $C_3$, by the map
$\varphi_3$, the corresponding point on $E_3$ is
\[W''=(12Q(a^2+2ab-2b^2), -108a(a-b)b(Q+1)Q).\]
By the group law, we get
\[\begin{split}
[2]W''=&\bigg(-12Q(2a^2(a^2+2ab-2b^2)Q^2+(8a^4+4a^3b-a^2b^2-6ab^3+3b^4)Q\\
&+2a^2(a^2+2ab-2b^2))/((Q+1)^2a^2),\\
&-108Q(a^4b(a-b)Q^4-2a^2(2a^4+3a^3b-a^2b^2-4ab^3+2b^4)Q^3\\
&+2(a^2+ab-b^2)(4a^4+a^3b-2ab^3+b^4)Q^2-2a^2(2a^4+3a^3b-a^2b^2\\
&-4ab^3+2b^4)Q+a^4b(a-b))/((Q+1)^2a^2)\bigg).
\end{split}\]

By the same method of above, we have infinitely many
$\mathbb{Q}(Q)$-rational solutions $(T,p)$ of Equation (2.7).

The last one is
\begin{equation}\label{Eq26}
z(z+a)(z+b)=\frac{r(r+a)(r+b)}{s(s+a)(s+b)}.
 \end{equation}
Let $z=T,$ and put $t=T(T+a)(T+b)$. By Equation (2.7), there are
infinitely many $\mathbb{Q}(t)$-rational solutions $(r,s)$ of
Equation (2.8). This completes the proof of Theorem 1.4.
\end{proof}

\textbf{Example 3.} The point $[2]W$ on $E_1$ leads to the solutions
of Equation (2.5):
\[\begin{split}
(x,y)=&\bigg((9T^6+12(a+b)T^5+(4a^2+17ab+4b^2)T^4+9ab(a+b)T^3\\
&+(2a^3b+7a^2b^2+2ab^3)T^2+2a^2b^2(a+b)T+a^3b^3)\\
&/(T(9T^4+12(a+b)T^3+(4a^2+17ab+4b^2)T^2+6ab(a+b)T+3a^2b^2)),\\
&((3T^2+2Ta+2Tb+ab)(-(3a+3b)T^3-2(a+b)^2T^2-2ab(a+b)T-a^2b^2)\\
&/(T(9T^4+12(a+b)T^3+(4a^2+17ab+4b^2)T^2+6ab(a+b)T+3a^2b^2)\bigg).
\end{split}\]

The point $[2]W'$ on $E_2$ leads to the solutions of Equation (2.6):
\[\begin{split}
(u,v)=&\bigg((-27(a+b)T^5-54(a+b)^2T^4-9(a+b)(4a^2+11ab+4b^2)T^3\\
&-(8a^4+68a^3b+93a^2b^2+68ab^3+8b^4)T^2-6ab(a+b)(2a^2+ab+2b^2)T\\
&+9a^3b^3)/(3(3T+2a+2b)(9T^4+12(a+b)T^3+(4a^2+17ab+4b^2)T^2\\
&+6ab(a+b)T+3a^2b^2)),\\
&-(81T^6+216(a+b)T^5+(216a^2+513ab+216b^2)T^4\\
&+3(a+b)(32a^2+109ab+32b^2)T^3+(16a^4+136a^3b+267a^2b^2\\
&+136ab^3+16b^4)T^2+6ab(a+b)(2a^2+7ab+2b^2)T+9a^3b^3)\\
&/(3(3T+2a+2b)(9T^4+12(a+b)T^3+(4a^2+17ab+4b^2)T^2\\
&+6ab(a+b)T+3a^2b^2))\bigg).
\end{split}\]

The point $[2]W''$ on $E_3$ leads to the solutions of Equation
(2.7):
\[\begin{split}
(T,p)=\bigg(\frac{a^2b(-Qa+Qb+a)(Q+1)}{-a^3Q^2+(2a^3-3ab^2+b^3)Q-a^3},\frac{-a^2b(-Qa+a-b)(Q+1)}{-a^3Q^2+(2a^3-3ab^2+b^3)Q-a^3}\bigg),
\end{split}\]
where $Q=q(q+a)(q+b)$.

By Equation (2.7), we can get a rational parametric solution of
Equation (2.8):
\[(r,s)=\bigg(\frac{a^2b(-ta+tb+a)(t+1)}{-a^3t^2+(2a^3-3ab^2+b^3)t-a^3},\frac{-a^2b(-ta+a-b)(t+1)}{-a^3t^2+(2a^3-3ab^2+b^3)t-a^3}\bigg),\]
where $t=T(T+a)(T+b)$.

Let $a=1,b=3,$ combining the above solutions, we can get a rational
parametric solution of Equation (1.2):
\[\begin{split}
&z=T=-\frac{3(2Q+1)(Q+1)}{(Q-1)^2},\\
&x=\frac{9T^6+48T^5+91T^4+108T^3+123T^2+72T+27}{T(9T^4+48T^3+91T^2+72T+27)},\\
&y=-\frac{(3T^2+8T+3)(12T^3+32T^2+24T+9)}{T(9T^4+48T^3+91T^2+72T+27)},\\
&u=\frac{108T^5+864T^4+2628T^3+3533T^2+1656T-243}{3(3T+8)(9T^4+48T^3+91T^2+72T+27)},\\
&v=\frac{81T^6+864T^5+3699T^4+7764T^3+7795T^2+2952T+243}{-3(3T+8)(9T^4+48T^3+91T^2+72T+27)},\\
&p=-\frac{3(Q+2)(Q+1)}{(Q-1)^2},\\
&q=q,\\
&r=-\frac{3(2T^3+8T^2+6T+1)(T^3+4T^2+3T+1)}{(T^3+4T^2+3T-1)^2},\\
&s=-\frac{3(T^3+4T^2+3T+2)(T^3+4T^2+3T+1)}{(T^3+4T^2+3T-1)^2},
\end{split}\]
where $Q=q(q+1)(q+3)$, and $q$ is a rational parameter.

\section{Some related questions}

In 1986, S. Hirose \cite{Hirose} conjectured that for $n\neq4$ if
$(n-2)P^n_k-(n-4)=2P^n_l$, then $P^n_{P^n_k}=P^n_kP^n_l$ can be
expressed as the sum and difference of two other $n$-gonal numbers.
It is difficult to prove it. Following this idea, for $n=12,$ we
find an example:
\[\begin{split}
P^{12}_{215666848}
&=P^{12}_{33841736}+P^{12}_{212995132}\\
&=P^{12}_{2907011822107606}-P^{12}_{2907011822107598}\\
&=P^{12}_{6568}P^{12}_{14686}.
\end{split}\]

For general $n$-gonal numbers, we have an open question.

\begin{question}
Are there infinitely many $n$-gonal numbers, except $n=3,4,5,6,8$,
which at the same time can be written as the sum, difference, and
product of other $n$-gonal numbers?
\end{question}

In Theorem 1.1, we give infinitely many quadratic polynomials
$f(X)\in Q[X]$ such that Equation (1.1) has infinitely many integer
solutions, but it seems difficult to solve the following question.

\begin{question} Does there exist a polynomial $f(X)\in Q[X]$ with $\deg{f}\geq3$ such that Equation (1.1) or Equation (1.2)
has infinitely many integer solutions?
\end{question}

For polynomials $f(X)\in Q[X]$ with $\deg{f}\geq4$, we raise

\begin{question}
Does there exist a polynomial $f(X)\in Q[X]$ with $\deg{f}\geq4$
such that Equation (1.1) or Equation (1.2) has a nontrivial rational
solution?
\end{question}

\end{document}